\documentclass[10pt]{amsart}
\usepackage{latexsym}   
\usepackage{amssymb}    
\usepackage{amsmath}    
\usepackage{amsthm}     
\usepackage[all]{xy}
\usepackage{hyperref}
\newcommand{\pa}{\partial}

\newcommand{\del}{\delta}
\newcommand{\Del}{\Delta}\newcommand{\ep}{\epsilon}

\newcommand{\na}{\nabla}
\newcommand{\ti}{\tilde}

\renewcommand{\thefootnote}

\newtheorem{theorem}{Theorem}[section]

\newtheorem{proposition}[theorem]{Proposition}

\theoremstyle{definition}

\theoremstyle{remark}
\newtheorem{remark}[theorem]{Remark}

\numberwithin{equation}{section}

\title[Thread configurations for ellipsoids] {
Thread configurations for ellipsoids}

\author[  Ion I. Dinc\u{a}]{Ion I. Dinc\u{a}}
\address{Faculty of Mathematics and Informatics,
University of Bucharest,  14 Academiei Str., 010014, Bucharest,
Romania}
 \email{dinca@gta.math.unibuc.ro}
\thanks{Supported by the University of Bucharest}
\subjclass[2000]{}

\begin{document}

\keywords{(confocal) quadrics, thread construction}

\begin{abstract}
We discuss Darboux-Staude type of thread configurations for the
ellipsoid similar to Chasles-Graves type of thread configurations
for the ellipse. These threads are formed by rectilinear segments,
geodesic and line of curvature segments on the considered
ellipsoid and with tangents tangent to the given ellipsoid and a
fixed confocal hyperboloid with one sheet and preserve constant
length when the vertices of the configuration move on confocal
ellipsoids.
\end{abstract}

\maketitle

\tableofcontents \pagenumbering{arabic}

\section{Introduction}

In trying to familiarize ourselves with classical results about
confocal quadrics, Staude's 1882 thread construction of confocal
ellipsoids appears as an important result. In Hilbert-Cohn-Vossen
(\cite{HCV},\S 4) this construction appears as a generalization of
thread construction of confocal ellipses, the two foci (singular
sets of singular ellipses and hyperbolas of the confocal family
and by means of which one jumps from one type of confocal conics
to the other) being replaced with orthogonal focal curves (an
ellipse and a hyperbola in orthogonal planes and with vertices of
the one being the foci of the other) and thus singular sets of
singular ellipsoids (hyperboloids) of the confocal family and by
means of which one jumps from one type of confocal quadrics to
another.

Note however that Staude's original construction as it appears in
Salmon (\cite{S1},\S 421a,b) allows the initial fixed ellipsoid
and hyperboloid with one sheet of the confocal family to be
nonsingular, the singular case being obtained as a limit (also in
the last sections \S 397-\S 421 of Salmon \cite{S1} other type of
thread configurations for quadrics are discussed; for example
Chasles's result about a thread fixed at two points on a quadric,
stretched with a pen and the pen thus moving on a line of
curvature of a quadric confocal to the given one).

According to Salmon (\cite{S1},\S 421a,b) this non-singular thread
construction of confocal ellipsoids can be interpreted as the
generalization of Graves's non-singular thread construction of
confocal ellipses using a thread passed around a given ellipse.

According to the Editor Reginald A. P. Rogers' preface to Salmon
\cite{S1} (referring to Staude's thread construction of confocal
quadrics) {\it 'In the Golden Age of Euclidean geometry, analogues
of these types were of great interest to men like Jacobi,
MacCullagh, Chasles and M. Roberts, but Staude's constructions
have virtually brought the subject to a conclusion. Staude's
treatment is also an excellent illustration of the elementary and
visible meaning of elliptic and hyper-elliptic integrals.'}

Unfortunately Staude's original paper \cite{S2} on addition of
hyper-elliptic integrals is beyond our grasp, mainly due to our
ignorance of German; however Darboux's generalization from 1870 of
Chasles's result from moving polygons circumscribed to an ellipse,
with vertices situated on confocal ellipses (and thus of constant
perimeter) to moving polygons circumscribed to a quadric, with
vertices situated on other quadrics confocal to the given one (and
thus of constant perimeter) seems to be discussed at length there,
so it may be the case that Staude himself has other results about
thread configurations for quadrics (see also footnote in Darboux
(\cite{D1}, Vol {\bf 2},\S 466)); in fact it is the similarities
of the arguments from Salmon (\cite{S1},\S 421 a,b) with those of
Darboux (\cite{D1},Vol{\bf 2}, Livre IV,Ch XIV) that drew our
attention to this project.

In Coolidge (\cite{C}, Ch. XIII) Staude's thread construction of
confocal quadrics is generalized to non-Euclidean geometries
(space forms).

For Chasles's result the basic result is a theorem due to Graves
(an identity involving an elliptic integral which at the geometric
level boils down to the excess between the sum of the lengths of
tangents to the given ellipse from a point on an ellipse confocal
to the given one and the subtended arc of the given ellipse being
independent of the point on the confocal ellipse); conversely
Chasles's result implies Graves's.

For Darboux's generalization of Chasles's result to dimension
three the differential equation of a line tangent to two confocal
quadrics and of its linear element in elliptic coordinates plays a
fundamental r\^{o}le. On one hand the linear element is a perfect
square, which allows separation and separate accounting of the
elliptic variables and on the other hand straight lines, having
tangents tangent to two confocal quadrics (Chasles-Jacobi), have
linear element amenable to these type of computations; moreover
reflections in confocal quadrics (which appear at the vertices of
thread configurations) are accounted just by changing the sign of
the variation of the corresponding elliptic coordinate; thus
threads formed by rectilinear segments tangent to the given
ellipsoid and hyperboloid with one sheet will preserve constant
length when the vertices will move on ellipsoids confocal to the
given one. However Darboux's result is not of a general nature; it
requires certain rationality conditions similar to the rationality
conditions required by closed geodesics on ellipsoids.

Note that the same theorem of Chasles-Jacobi allows the
computations of lines in elliptic coordinates to be extended to
geodesic segments with tangents tangent to the same two confocal
quadrics; their part involving the linear element is also
trivially extended to segments of intersections of the two given
quadrics (lines of curvature). While line of curvature segments
are not locally length minimizing under the condition of being
situated on one side of a quadric, they are situated on the same
side of two quadrics and thus are allowed in thread configurations
as boundary requirements; if one (or both) quadric(s) becomes
singular, then the line of curvature segments become geodesic
segments and we have a genuine variational problem.

Staude's thread construction of the ellipsoid as it appears in
Salmon (\cite{S1},\S 421a,b) has only one vertex and two (either
possibly void) line of curvature segments; thus in this vein
Chasles's and Graves's result are in a relation analogous to that
of Darboux's and Staude's.

We can extend threads to allow them to be formed by rectilinear
segments, geodesic segments and line of curvature segments with
common tangent at the points of change from one type of segment to
the other (thus the threads will be analytic on pieces and
(excepting vertices) with continuous derivative); keeping in mind
variations of elliptic coordinates (certain of their extreme
values are alternatively attained) their length will be constant.

For threads without line of curvature segments a rationality
condition remains; if we allow line of curvature segments, then
the rationality condition disappears.

We also derive the algebraic structure of vertex configuration via
the Ivory affinity between confocal quadrics.

Using this local result and Darboux's constant perimeter property
concerning moving polygons circumscribed to a given set of $n$ and
inscribed in arbitrarily many $n$-dimensional confocal quadrics we
derive Darboux's result as a variational principle, valid in a
more general form (complex settting) and which covers all (totally
real) quadrics in the complex Euclidean space.

A simple internet search with keywords reveals most of the current
literature in this area.

There has been some work of J. Itoh and K. Kiyohara on simplifying
and generalizing Staude's thread construction of confocal quadrics
to higher dimensional quadrics and space forms as ambient space,
but unfortunately we have been unable to get free access to it.

There is also some recent work on Chasles-Darboux type results,
related to closed billiard trajectories, confocal quadrics and
hyper-elliptic integrals (see Dragovi\'{c}-Radnovi\'{c} \cite{DR2}
for a synthetic approach, Tabachnikov \cite{T} and their
references) and to closed geodesics on the ellipsoid (classically
studied by Jacobi, Weierstrass, etc and recently by Knorrer,
Moser, etc; see Abenda-Fedorov \cite{AF}, Fedorov \cite{F} and
their references). For example in Darboux's result when the
vertices are situated all on the same ellipsoid the polygon in
question can be viewed as a closed geodesic on a degenerate
$3$-dimensional ellipsoid (the double cover of the interior of the
$2$-dimensional ellipsoid in question) and closed geodesics on
ellipsoids require certain rationality conditions; these
rationality conditions remain valid in the degenerate case and are
precisely those found by Darboux.

\section{Confocal quadrics in canonical form}

Consider the complexified Euclidean space
$$(\mathbb{C}^{n+1},<.,.>),\
<x,y>:=x^Ty,\ |x|^2:=x^Tx,\ x,y\in\mathbb{C}^{n+1}$$ with standard
basis $\{e_j\}_{j=1,...,n+1},\ e_j^Te_k=\del_{jk}$.

Isotropic (null) vectors are those vectors $v$ of length $0\
(|v|^2=0)$; since most vectors are not isotropic we shall call a
vector simply vector and we shall only emphasize isotropic when
the vector is assumed to be isotropic. The same denomination will
apply in other settings: for example we call quadric a
non-degenerate quadric (a quadric projectively equivalent to the
complex unit sphere).

A quadric $x\subset\mathbb{C}^{n+1}$ is given by the quadratic
equation
$$Q(x):=\begin{bmatrix}x\\1\end{bmatrix}^T\begin{bmatrix}A&B\\B^T&C\end{bmatrix}
\begin{bmatrix}x\\1\end{bmatrix}=x^T(Ax+2B)+C=0,$$
$$A=A^T\in\mathbf{M}_{n+1}(\mathbb{C}),\
B\in\mathbb{C}^{n+1},\ C\in\mathbb{C},\
\begin{vmatrix}A&B\\B^T&C\end{vmatrix}\neq 0.$$

There are many definitions of totally real (sub)spaces of
$\mathbb{C}^{n+1}$, some even involving a hermitian inner product,
but all definitions coincide: an $(n+1)$-totally real subspace of
$\mathbb{C}^{n+1}$ is of the form
$(R,t)(\mathbb{R}^k\times(i\mathbb{R})^{n+1-k}),\ k=0,...,n+1$,
where
$(R,t)\in\mathbf{O}_{n+1}(\mathbb{C})\ltimes\mathbb{C}^{n+1}$. Now
a totally real quadric is simply an $n$-dimensional quadric in an
$(n+1)$-totally real subspace of $\mathbb{C}^{n+1}$.

A metric classification of all (totally real) quadrics in
$\mathbb{C}^{n+1}$ requires the notion of {\it symmetric Jordan}
(SJ) canonical form of a symmetric complex matrix. The symmetric
Jordan blocks are:
$$J_1:=0=0_{1,1}\in\mathbf{M}_1(\mathbb{C}),\
J_2:=f_1f_1^T\in\mathbf{M}_2(\mathbb{C}),\
J_3:=f_1e_3^T+e_3f_1^T\in\mathbf{M}_3(\mathbb{C}),$$
$$J_4:=f_1\bar
f_2^T+f_2f_2^T+\bar f_2f_1^T\in\mathbf{M}_4(\mathbb{C}),\ J_5:=
f_1\bar f_2^T+f_2e_5^T+e_5f_2^T+\bar
f_2f_1^T\in\mathbf{M}_5(\mathbb{C}),$$
$$J_6:= f_1\bar f_2^T+f_2\bar
f_3^T+f_3f_3^T+\bar f_3f_2^T+\bar
f_2f_1^T\in\mathbf{M}_6(\mathbb{C}),$$ etc, where
$f_j:=\frac{e_{2j-1}+ie_{2j}}{\sqrt{2}}$ are the standard
isotropic vectors (at least the blocks $J_2,\ J_3$ were known to
the classical geometers). Any symmetric complex matrix can be
brought via conjugation with a complex rotation to the symmetric
Jordan canonical form, that is a matrix block decomposition with
blocks of the form $a_jI_p+J_p$; totally real quadrics are
obtained for eigenvalues $a_j$ of the quadratic part $A$ defining
the quadric being real or coming in complex conjugate pairs $a_j,\
\bar a_j$ with subjacent symmetric Jordan blocks of same dimension
$p$. Just as the usual Jordan block $\sum_{j=1}^pe_je_{j+1}^T$ is
nilpotent with $e_{p+1}$ cyclic vector of order $p$, $J_p$ is
nilpotent with $\bar f_1$ cyclic vector of order $p$, so we can
take square roots of SJ matrices without isotropic kernels
($\sqrt{aI_p+J_p}:=\sqrt{a}\sum_{j=0}^{p-1}(^{\frac{1}{2}}_j)a^{-j}J_p^j,\
a\in\mathbb{C}^*,\ \sqrt{a}:=\sqrt{r}e^{i\theta}$ for
$a=re^{2i\theta},\ 0<r,\ -\pi<2\theta\le\pi$), two matrices with
same SJ decomposition type (that is $J_p$ is replaced with a
polynomial in $J_p$) commute, etc.

The confocal family $\{x_z\}_{z\in\mathbb{C}}$ of a quadric
$x_0\subset\mathbb{C}^{n+1}$ in canonical form (depending on as
few constants as possible) is given in the projective space
$\mathbb{C}\mathbb{P}^{n+1}$ by the equation
$$Q_z(x_z):=\begin{bmatrix}x_z\\1\end{bmatrix}^T
(\begin{bmatrix}A&B\\B^T&C\end{bmatrix}^{-1}-z
\begin{bmatrix}I_{n+1}&0\\0^T&0\end{bmatrix})^{-1}\begin{bmatrix}x_z\\1\end{bmatrix}=0,$$
where

$\bullet\ \  A=A^T\in\mathbf{GL}_{n+1}(\mathbb{C})$ SJ,
$B=0\in\mathbb{C}^{n+1},\ C=-1$ for {\it quadrics with center}
(QC),

$\bullet\ \ A=A^T\in\mathbf{M}_{n+1}(\mathbb{C})$ SJ,
$\ker(A)=\mathbb{C}e_{n+1},\ B=-e_{n+1},\ C=0$ for {\it quadrics
without center} (QWC) and

$\bullet\ \ A=A^T\in\mathbf{M}_{n+1}(\mathbb{C})$ SJ,
$\ker(A)=\mathbb{C}f_1,\ B=-\bar f_1,\ C=0$ for {\it isotropic
quadrics without center} (IQWC).

The totally real confocal family of a totally real quadric is
obtained for $z\in\mathbb{R}$.

From the definition one can see that the family of quadrics
confocal to $x_0$ is the adjugate of the pencil generated by the
adjugate of $x_0$ and Cayley's absolute
$C(\infty)\subset\mathbb{C}\mathbb{P}^n$ in the hyperplane at
infinity; since Cayley's absolute encodes the Euclidean structure
of $\mathbb{C}^{n+1}$ (it is the set invariant under rigid motions
and homotheties of
$\mathbb{C}^{n+1}:=\mathbb{CP}^{n+1}\backslash\mathbb{CP}^n$) the
mixed {\it metric-projective} character of the confocal family
becomes clear.

For QC $\mathrm{spec}(A)$ is unambiguous (does not change under
rigid motions
$(R,t)\in\mathbf{O}_{n+1}(\mathbb{C})\ltimes\mathbb{C}^{n+1}$) but
for (I)QWC it may change with $(p+1)$-roots of unity for the block
of ($f_1$ in $A$ being $J_p$) $e_{n+1}$ in $A$ being $J_1$ even
under rigid motions which preserve the canonical form, so it is
unambiguous up to $(p+1)$-roots of unity; for simplicity we make a
choice and work with it.

We have the diagonal Q(W)C respectively for
$A=\Sigma_{j=1}^{n+1}a_j^{-1}e_je_j^T,\
A=\Sigma_{j=1}^na_j^{-1}e_je_j^T$; the diagonal IQWC come in
different flavors, according to the block of $f_1:\
A=J_p+\Sigma_{j=p+1}^{n+1}a_j^{-1}e_je_j^T$; in particular if
$A=J_{n+1}$, then $\mathrm{spec}(A)=\{0\}$ is unambiguous. General
quadrics are those for which all eigenvalues have geometric
multiplicity $1$; equivalently each eigenvalue has an only
corresponding SJ block; in this case the quadric also admits
elliptic coordinates.

There are continuous groups of symmetries which preserve the SJ
canonical form for more than one SJ block corresponding to an
eigenvalue, so from a metric point of view a metric classification
according to the elliptic coordinates and continuous symmetries
may be a better one.

With $R_z:=I_{n+1}-zA,\
z\in\mathbb{C}\setminus\mathrm{spec}(A)^{-1}$ the family of
quadrics $\{x_z\}_z$ confocal to $x_0$ is given by
$Q_z(x_z)=x_z^TAR_z^{-1}x_z+2(R_z^{-1}B)^Tx_z+C+zB^TR_z^{-1}B=0$.
For $z\in\mathrm{spec}(A)^{-1}$ we obtain singular confocal
quadrics; those with $z^{-1}$ having geometric multiplicity $1$
admit a singular set which is an $(n-1)$-dimensional quadric
projectively equivalent to $C(\infty)$, so they will play an
important r\^{o}le in the discussion of homographies
$H\in\mathbf{PGL}_{n+1}(\mathbb{C})$ taking a confocal family into
another one, since $H^{-1}(C(\infty)),\ C(\infty)$ respectively
$C(\infty),\ H(C(\infty))$ will suffice to determine each confocal
family. Such homographies preserve all metric-projective
properties of confocal quadrics (including the good metric
properties of the Ivory affinity) and thus all integrable systems
whose integrability depends only on the family of confocal
quadrics (the resulting involutory transformation between
integrable systems is called {\it Hazzidakis} (H) by Bianchi for
the integrable system in discussion being the problem of deforming
quadrics). While the spectrum ($z$'s of $Q_z$) of a family of
confocal quadrics is not well defined, the relative spectrum
(difference of $z$'s) is; thus we can consider $Q_z$ as $Q_0$ for
any $z\in\mathbb{C}$ by a translation of $\mathbb{C}$ which brings
$z$ to $0$.

\subsection{Some classical metric properties of confocal quadrics and of
the Ivory affinity}

The Ivory affinity is an affine correspondence between confocal
quadrics and having good metric properties (it may be the reason
why Bianchi calls it {\it affinity} in more than one language): it
is given by $x_z=\sqrt{R_z}x_0+C(z),\
C(z):=-(\frac{1}{2}\int_0^z(\sqrt{R_w})^{-1}dw)B$. Note that
$C(z)=0$ for QC, $=\frac{z}{2}e_{n+1}$ for QWC; for IQWC it is the
Taylor series of $\frac{1}{2}\int_0^z(\sqrt{1-w})^{-1}dw$ at $z=0$
with each monomial $z^{k+1}$ replaced by $z^{k+1}J_p^k\bar f_1$,
where $J_p$ is the block of $f_1$ in $A$ and thus a polynomial of
degree $p$ in $z$. Note
$AC(z)+(I_{n+1}-\sqrt{R_z})B=0=(I_{n+1}+\sqrt{R_z})C(z)+zB$ (both
are $0$ for $z=0$ and do not depend on $z$). Applying $d$ to
$Q_z(x_z)=0$ we get $dx_z^TR_z^{-1}(Ax_z+B)=0$, so the unit normal
$N_z$ is proportional to $\hat N_z:=-2\pa_zx_z$. If
$\mathbb{C}^{n+1}\ni x\in x_{z_1},x_{z_2}$, then $\hat
N_{z_j}=R_{z_j}^{-1}(Ax+B)$; using $R_z^{-1}-I_{n+1}=zAR_z^{-1},\
z_1R_{z_1}^{-1}-z_2R_{z_2}^{-1}=(z_1-z_2)R_{z_1}^{-1}R_{z_2}^{-1}$
we get $0=Q_{z_1}(x)-Q_{z_2}(x)=(z_1-z_2)\hat N_{z_1}^T\hat
N_{z_2}$, so two confocal quadrics cut each other orthogonally
(Lam\'{e}). For general quadrics the polynomial equation
$Q_z(x)=0$ has degree $n+1$ in $z$ and it has multiple roots {\it
if and only if} (iff) $0=\pa_zQ_z(x)=|\hat N_z|^2$; thus outside
the locus of isotropic normals elliptic coordinates (given by the
roots $z_1,...,z_{n+1}$ of the said equation) give a
parametrization of $\mathbb{C}^{n+1}$ suited to confocal quadrics.
With $x_0^0,x_0^1\in x_0,\ V_0^1:=x_z^1-x_0^0$, etc the
preservation of length of segments between confocal quadrics
(Ivory Theorem) becomes
$|V_0^1|^2=|x_0^0+x_0^1-C(z)|^2-2(x_0^0)^T(I_{n+1}+\sqrt{R_z})x_0^1+zC=|V_1^0|^2$;
the preservation of lengths of rulings (Henrici):
$w_0^TAw_0=w_0^T\hat N_0=0,\ w_z=\sqrt{R_z}w_0\Rightarrow
w_z^Tw_z=|w_0|^2-zw_0^TAw_0=|w_0|^2$; the symmetry of the {\it
tangency configuration} (TC)(Bianchi): $(V_0^1)^T\hat N_0^0
=(x_0^0)^TA\sqrt{R_z}x_0^1-B^T(x_z^0+x_z^1-C(z))+C=(V_1^0)^T\hat
N_0^1$; the preservation of angles between segments and rulings
(Bianchi): $(V_0^1)^Tw_0^0+(V_1^0)^Tw_z^0=-z(\hat
N_0^0)^Tw_0^0=0$; the preservation of angles between rulings
(Bianchi):
$(w_0^0)^Tw_z^1=(w_0^0)^T\sqrt{R_z}w_0^1=(w_z^0)^Tw_0^1$; the
preservation of angles between polar rulings: $(w_0^0)^TA\hat
w_0^0=0\Rightarrow (w_z^0)^T\hat w_z^0=(w_0^0)^T\hat
w_0^0-z(w_0^0)^TA\hat w_0^0=(w_0^0)^T\hat w_0^0$, etc.

We also have the Chasles-Jacobi result (that of Jacobi's being an
inspiration for that of Chasles's)

{\it (Jacobi) The tangent lines to a geodesic on $x_0$ remain
tangent to $n-1$ other confocal quadrics.}

{\it (Chasles) The common tangents to $n$ confocal quadrics form a
normal congruence and envelope geodesics on the $n$ confocal
quadrics.}

\section{Chasles's and Graves's results}

Let $a_1>a_2>0$ and consider the confocal ellipses
$x_z:=[\sqrt{a_1-z}\cos\theta\ \ \sqrt{a_2-z}\sin\theta]^T,\ z<0,\
\theta\in\mathbb{R}$ outside the initial ellipse $x_0$.

Chasles's result about polygons circumscribed to a given ellipse
and inscribed in a given set of ellipses confocal to the given one
roughly states:

\begin{theorem} (Chasles)

Given a set of ellipses $x_{z_j},\ z_j<0,\ j=1,...,n$ confocal to
the given one $x_0$, if a ray of light tangent to $x_0$ at a point
$x_0^1\in x_0$ be reflected in $x_{z_1}$, then it remains tangent
to $x_0$ at $x_0^2\in x_0$ and if after successive reflections
further in $x_{z_2},...,x_{z_n}$ it returns to being tangent to
$x_0$ at $x_0^1$, then this property and the perimeter of the
obtained polygon is independent of the position of $x_0^1$ on
$x_0$.

\end{theorem}

In fact Chasles proved also a dual result related to a ray of
light reflecting on $x_0$ and tangent (possibly outside the
segment of the actual trajectory) to confocal hyperbolas
$x_{z_j},\ a_1>z_j>a_2,\ j=1,...,n$; should one of the $x_{z_j}$
be an ellipse with $a_2>z_j>0$ then the ray of light becomes
entrapped and will always be tangent to $x_{z_j}$ (the Poncelet
Theorem).

The fact that Chasles's result is an immediate consequence of the
next

\begin{theorem} (Graves)

If a thread longer than the perimeter of an ellipse $x_0$ is
passed around $x_0$ and stretched with a pen, then the tip of the
pen will describe an ellipse $x_z,\ z<0$ confocal with the given
one $x_0$

\end{theorem}

is straightforward, Chasles's result being a particular
configuration of $n$ threads (or equivalently of a single thread
stretched with $n$ pens) such that their rectilinear parts are in
continuation one of the other.

Conversely, Chasles's result implies the previous by taking
$z_1:=z<0$ fixed, $z_j<0,\ j=2,...,n$ infinitesimally close to $0$
and letting $n\rightarrow\infty$.

If we give up the requirement that the rectilinear parts are in
continuation one of the other, then we get more general thread
configurations of the ellipse: a closed thread around $x_0$ and
stretched with $n$ pens whose tips are situated on ellipses
$x_{z_1},...,x_{z_n}$ such that consecutive rectilinear parts are
in continuation one of the other or are joined by an arc of the
ellipse $x_0$; we may even allow an ideal construction that
consecutive rectilinear parts may not be in continuation one of
the other and are joined by an arc of the ellipse but with cusps
at both ends, in which case the arc of the ellipse joining their
ends at $x_0$ has to be subtracted to get constant perimeter.

The Ivory affinity between confocal ellipses is an affine
correspondence given by $x_z=\sqrt{R_z}x_0,\\ R_z:=I_2-zA,\
A:=\mathrm{diag}[a_1^{-1}\ \ a_2^{-1}]$; as $z$ varies $x_z$
describes orthogonal curves of the family of confocal ellipses
(the hyperbolas of the confocal family).

It has good metric properties, among which are the preservation of
lengths of segments between confocal ellipses (if $x_0^0,\
x_0^1\in x_0$ and by use of the Ivory affinity we get $x_z^0,\
x_z^1\in x_z$, then $|x_z^1-x_0^0|=|x_z^0-x_0^1|$) and of the {\it
tangency configuration} (TC) ($x_z^1-x_0^0$ is tangent to $x_0$ at
$x_0^0$ {\it if and only if} (iff) $x_z^0-x_0^1$ is tangent to
$x_0$ at $x_0^1$).

We have thus reduced our investigation to:

\begin{proposition} (Ivory affinity approach for Graves's vertex
configuration)

Given $z<0$ and three points $x_0^1,x_0^0,x_0^2\in x_0$ in this
order and by the Ivory affinity the corresponding
$x_z^1,x_z^0,x_z^2\in x_z$, then $x_z^0-x_0^1,\ x_z^0-x_0^2$
reflect in $x_z$ at $x_z^0$ iff $x_z^1-x_0^0,\ x_z^2-x_0^0$
reflect in $x_0$ at $x_0^0$; further in this case, by the
preservation of the TC under the Ivory affinity, $x_z^0-x_0^1$ is
tangent to $x_0$ at $x_0^1$ iff $x_z^0-x_0^2$ is tangent to $x_0$
at $x_0^2$ and further in this case the excess
$|x_z^1-x_z^2|-\mathrm{length}_{x_0}(x_0^1x_0^0x_0^2)$ does not
depend on the position of $x_0^0$ on $x_0$.

\end{proposition}

\begin{proof}

The condition that $x_0^1-x_z^0,\ x_0^2-x_z^0$ reflect in $x_z$ at
$x_z^0$ becomes
$0=(dx_z^0)^T(\frac{x_0^1-x_z^0}{|x_0^1-x_z^0|}+\frac{x_0^2-x_z^0}{|x_0^2-x_z^0|})=
(dx_0^0)^T(\frac{x_z^1-x_0^0}{|x_z^1-x_0^0|}+\frac{x_z^2-x_0^0}{|x_z^2-x_0^0|})$;
the TC follows because the segment $[x_z^1,x_z^2]$ must be tangent
to $x_0$ at $x_0^0$. Simple variational arguments show that the
reflection property is enough to infer the constant excess
property in the TC (the gradient lines of the thread length
function are the confocal hyperbolas of the confocal family);
however we shall explicitly state the constant excess property as
an identity involving an elliptic integral and show how the
reflection property is used at the analytic level also.

With $x_0^j=[\sqrt{a_1}\cos\theta_j\ \ \sqrt{a_2}\sin\theta_j],\
j=0,1,2$ we have the TC
\begin{eqnarray}\label{eq:tc}
\cos\theta_0\sqrt{1-za_1^{-1}}\cos\theta_{1,2}+\sin\theta_0\sqrt{1-za_2^{-1}}\sin\theta_{1,2}=1,
\end{eqnarray} so
\begin{eqnarray}\label{eq:t12}
\cos\theta_j=\frac{\cos\theta_0\sqrt{1-za_1^{-1}}-(-1)^j\sin\theta_0\sqrt{1-za_2^{-1}}
\sqrt{-z(a_2^{-1}\sin^2\theta_0+a_1^{-1}\cos^2\theta_0)}}
{1-z(a_2^{-1}\sin^2\theta_0+a_1^{-1}\cos^2\theta_0)},\nonumber\\
\sin\theta_j=\frac{\sin\theta_0\sqrt{1-za_2^{-1}}+(-1)^j\cos\theta_0\sqrt{1-za_1^{-1}}
\sqrt{-z(a_2^{-1}\sin^2\theta_0+a_1^{-1}\cos^2\theta_0)}}
{1-z(a_2^{-1}\sin^2\theta_0+a_1^{-1}\cos^2\theta_0)},\nonumber\\
j=1,2.
\end{eqnarray}
We now need
$0=\frac{d}{d\theta_0}[\frac{|x_z^1-x_z^2|}{\sqrt{a_1a_2}}-
\int_{\theta_1}^{\theta_2}\sqrt{a_2^{-1}\sin^2\theta+a_1^{-1}\cos^2\theta}d\theta]
=\frac{1}{z}\frac{d}{d\theta_0}\frac{2\sqrt{-z(1-za_1^{-1})(1-za_2^{-1})}}
{1-z(a_2^{-1}\sin^2\theta_0+a_1^{-1}\cos^2\theta_0)}\\-
\sqrt{a_2^{-1}\sin^2\theta_2+a_1^{-1}\cos^2\theta_2}\frac{d\theta_2}{d\theta_0}
+\sqrt{a_2^{-1}\sin^2\theta_1+a_1^{-1}\cos^2\theta_1}\frac{d\theta_1}{d\theta_0}$,
or, using (\ref{eq:tc}) for $\frac{d\theta_j}{d\theta_0},\\
\sin\theta_0\cos\theta_0(a_2^{-1}-a_1^{-1})\sqrt{-z}
(\frac{4\sqrt{(1-za_1^{-1})(1-za_2^{-1})}}
{1-z(a_2^{-1}\sin^2\theta_0+a_1^{-1}\cos^2\theta_0)}
-\frac{\sum_{j=1}^2\sqrt{a_2^{-1}\sin^2\theta_j+a_1^{-1}\cos^2\theta_j}}
{\sqrt{a_2^{-1}\sin^2\theta_0+a_1^{-1}\cos^2\theta_0}})\\
-\sqrt{(1-za_1^{-1})(1-za_2^{-1})}
\sum_{j=1}^2(-1)^j\sqrt{a_2^{-1}\sin^2\theta_j+a_1^{-1}\cos^2\theta_j}=0$.

Now from the reflection property we have
$0=(\frac{dx_z^0}{d\theta_0})^T\sum_{j=1}^2\frac{x_0^j-x_z^0}{|x_0^j-x_z^0|}=\\
\sum_{j=1}^2\frac{(a_1-a_2)\sin\theta_0\cos\theta_0+a_2\sin\theta_j\cos\theta_0\sqrt{1-za_2^{-1}}
-a_1\cos\theta_j\sin\theta_0\sqrt{1-za_1^{-1}}}
{\sqrt{a_1(\cos\theta_j-\cos\theta_0\sqrt{1-za_1^{-1}})^2+
a_2(\sin\theta_j-\sin\theta_0\sqrt{1-za_2^{-1}})^2}}$, or, using
(\ref{eq:t12}):
$\\\sum_{j=1}^2\frac{\sin\theta_0\cos\theta_0(a_2^{-1}-a_1^{-1})\sqrt{-z}
+(-1)^j\sqrt{(a_2^{-1}\sin^2\theta_0+a_1^{-1}\cos^2\theta_0)(1-za_1^{-1})(1-za_2^{-1})}}
{\sqrt{a_2^{-1}\sin^2\theta_j+a_1^{-1}\cos^2\theta_j}}=0$, so we
need $\\\sqrt{a_2^{-1}\sin^2\theta_j+a_1^{-1}\cos^2\theta_j}
=\frac{\sqrt{(a_2^{-1}\sin^2\theta_0+a_1^{-1}\cos^2\theta_0)(1-za_1^{-1})(1-za_2^{-1})}
+(-1)^j\sin\theta_0\cos\theta_0(a_2^{-1}-a_1^{-1})\sqrt{-z}}
{1-z(a_2^{-1}\sin^2\theta_0+a_1^{-1}\cos^2\theta_0)},\\ j=1,2$,
which is straightforward.

\end{proof}

\section{Darboux's and Staude's results}

Here we shall reproduce computations mainly from Darboux
(\cite{D1},Vol{\bf 2}, Livre IV,Ch XIV) and Salmon (\cite{S1},\S
421a,b) concerning elliptic coordinates on $\mathbb{R}^3$ and
their corresponding applications (see also Bianchi (\cite{B},\S
419-\S 427)).

\subsection{Elliptic coordinates on $\mathbb{R}^3$}\noindent

Let $a_1>a_2>a_3$ constants and $a_1>u^1>a_2>u^2>a_3>u^3$ be the
elliptic coordinates on $\mathbb{R}^3\setminus\{x^1x^2x^3=0\}\ni
x=[x^1\ \ x^2\ \ x^3]^T,\
(x^j)^2=\frac{\prod_k(a_j-u^k)}{\prod_{k\neq j}(a_j-a_k)},\
j=1,2,3$.

Because both the numerator and the denominator in each term on the
right hand side contain the same number of negative terms, this
provides good definition.

Also we have
\begin{eqnarray}\label{eq:ell}
\sum_j\frac{(x^j)^2}{a_j-u}-1=\frac{\prod_j(u-u^j)}{\prod_j(a_j-u)},
\end{eqnarray}
namely this relation is separately linear in $u^k,\ k=1,2$, so it
is sufficient to verify it for two different values of an $u^k$;
of course we take $u^k:=a_k$ and $u^k:=a_{k+1}$, in which case we
have reduced what we had to prove from dimension $3$ to dimension
$2$ and we have a backward induction.

We have the linear element of
$\mathbb{R}^3\setminus\{x^1x^2x^3=0\}$ in elliptic coordinates
\begin{eqnarray}\label{eq:lell}
ds^2=|dx|^2=\sum_k\frac{\prod_{j\neq
k}(u^k-u^j)}{4\prod_j(a_j-u^k)}(du^k)^2.
\end{eqnarray}

To see this first we need $\sum_j\frac{a_j-u}{\prod_{k\neq
j}(a_j-a_k)}=0$ (which follows by letting $u:=a_1,a_2$) and then
differentiate (\ref{eq:ell}) with respect to $u$ and let $u:=u^k$.

From (\ref{eq:ell}) for $a_1>u^1>a_2$ constant we get the
hyperboloids with two sheets of the confocal family, for
$a_2>u^2>a_3$ constant we get the hyperboloids with one sheet of
the confocal family and for $a_3>u^3>-\infty$ constant we get the
ellipsoids of the confocal family.

Note also that the singular cases $u^k=a_k,\ a_{k+1}$ are allowed
and make sense, both at the level of elliptic coordinates and
linear element, by a limiting argument $u^k\nearrow a_k,\
u^k\searrow a_{k+1}$ (of course for each such singular value of
$u^k$ we have a similar $2$-dimensional discussion for the
remaining parameters).

For $u^1\nearrow a_1$ we get the doubly covered plane $\{x^1=0\}$
with elliptic coordinates $[\pm 0\\
\pm\sqrt{\frac{(a_2-u^2)(a_2-u^3)}{a_2-a_3}}\ \
\pm\sqrt{\frac{(a_3-u^2)(a_3-u^3)}{a_3-a_2}}]^T,\ a_2>u^2>a_3>u^3$
as the limit of a fattening hyperboloid with two sheets whose
sheets tend to the point $0$ and this forces them to elongate
themselves along the plane $\{x^1=0\}$.

For $u^1\searrow a_2$ we get the doubly covered convex region of
the hyperbola $\frac{(x^1)^2}{a_1-a_2}+\frac{(x^3)^2}{a_3-a_2}=1$
in the plane $\{x^2=0\}$ with elliptic coordinates
$[\pm\sqrt{\frac{(a_1-u^2)(a_1-u^3)}{a_1-a_3}}\ \ \pm 0\ \
\pm\sqrt{\frac{(a_3-u^2)(a_3-u^3)}{a_3-a_1}}]^T,\ a_2>u^2>a_3>u^3$
as the limit of a thinning hyperboloid with two sheets which tends
to the plane $\{x^2=0\}$.

For $u^2\nearrow a_2$ we get the doubly covered concave region of
the hyperbola $\frac{(x^1)^2}{a_1-a_2}+\frac{(x^3)^2}{a_3-a_2}=1$
in the plane $\{x^2=0\}$ with elliptic coordinates
$[\pm\sqrt{\frac{(a_1-u^1)(a_1-u^3)}{a_1-a_3}}\ \ \pm 0\ \
\pm\sqrt{\frac{(a_3-u^1)(a_3-u^3)}{a_3-a_1}}]^T,\ a_1>u^1>a_2,\
a_3>u^3$ as the limit of a thinning (without hole) hyperboloid
with one sheet which tends to the plane $\{x^2=0\}$; thus for
$u^1\searrow a_2\nwarrow u^2$ we get the focal hyperbola
$[\pm\sqrt{\frac{(a_1-a_2)(a_1-u^3)}{a_1-a_3}}\\ \pm 0\ \
\pm\sqrt{\frac{(a_3-a_2)(a_3-u^3)}{a_3-a_1}}]^T$ covered once from
each side.

For $u^2\searrow a_3$ we get the doubly covered concave region of
the ellipse $\frac{(x^1)^2}{a_1-a_3}+\frac{(x^2)^2}{a_2-a_3}=1$ in
the plane $\{x^3=0\}$ with elliptic coordinates
$[\pm\sqrt{\frac{(a_1-u^1)(a_1-u^3)}{a_1-a_2}}\ \
\pm\sqrt{\frac{(a_2-u^1)(a_2-u^3)}{a_2-a_1}}\ \ \pm 0]^T,\
a_1>u^1>a_2,\ a_3>u^3$ as the limit of a thinning (with hole)
hyperboloid with one sheet which tends to the plane $\{x^3=0\}$.

For $u^3\nearrow a_3$ we get the doubly covered convex region of
the ellipse $\frac{(x^1)^2}{a_1-a_3}+\frac{(x^2)^2}{a_2-a_3}=1$ in
the plane $\{x^3=0\}$ with elliptic coordinates
$[\pm\sqrt{\frac{(a_1-u^1)(a_1-u^2)}{a_1-a_2}}\ \
\pm\sqrt{\frac{(a_2-u^1)(a_2-u^2)}{a_2-a_1}}\ \ \pm 0]^T,\
a_1>u^1>a_2>u^2>a_3$ as the limit of a thinning ellipsoid which
tends to the plane $\{x^3=0\}$; thus for $u^2\searrow a_3\nwarrow
u^3$ we get the focal ellipse
$[\pm\sqrt{\frac{(a_1-u^1)(a_1-a_3)}{a_1-a_2}}\ \
\pm\sqrt{\frac{(a_2-u^1)(a_2-a_3)}{a_2-a_1}}\ \ 0]^T$ covered once
from each side.

For $u^3\searrow -\infty$ the ellipsoid tends to infinity and in
shape closer to a sphere.

Thus through each point of $\mathbb{R}^3$ pass $3$ quadrics of the
confocal family which (excluding the focal ellipse and hyperbola)
meet orthogonally (Lam\'{e}).

We now fix a hyperboloid with one sheet $\{u^2=u^2_0\}$ and an
ellipsoid $\{u^3=u^3_0\}$; we are further interested in the region
of the space outside both of them:
$a_1>u^1>a_2>u^2_0>u^2>a_3>u^3_0>u^3$.

We are interested in the {\it congruence} ($2$-dimensional family)
of lines tangent to both $\{u^2=u^2_0\}$ and $\{u^3=u^3_0\}$ and
in the linear elements of these lines in elliptic coordinates;
according to the Chasles-Jacobi result this congruence is {\it
normal} (the lines are the normals to an $1$-dimensional family of
surfaces) and its developables envelope geodesics on the two
quadrics.

These surfaces are given in elliptic coordinates by
\begin{eqnarray}\label{eq:Phi}
\Phi:=\frac{1}{2}\sum_k\ep_k\int\sqrt{\frac{(u^k-u^2_0)(u^k-u^3_0)}{\prod_j(a_j-u^k)}}du^k=
\mathrm{ct},\ \ep_k:=\pm 1,\ k=1,2,3.
\end{eqnarray}
We have $|\na\Phi|^2=\sum_k\frac{4\prod_j(a_j-u^k)}{\prod_{j\neq
k}(u^k-u^j)}(\frac{\pa\Phi}{\pa
u^k})^2=\sum_k\frac{(u^k-u^2_0)(u^k-u^3_0)}{\prod_{j\neq
k}(u^k-u^j)}=1$ (again we can let for example $u^2_0:=u^1,u^2$),
so these surfaces are indeed parallel.

Differentiating this with respect to $u^2_0,u^3_0$ we obtain
$(\na\Phi)^T\na\frac{\pa\Phi}{\pa
u^2_0}=(\na\Phi)^T\na\frac{\pa\Phi}{\pa u^3_0}=0$; keeping account
also of $(\na\frac{\pa\Phi}{\pa u^2_0})^T\na\frac{\pa\Phi}{\pa
u^3_0}=\sum_k\frac{4\prod_j(a_j-u^k)}{\prod_{j\neq
k}(u^k-u^j)}\frac{\pa^2\Phi}{\pa u^k\pa u^2_0}\frac{\pa^2\Phi}{\pa
u^k\pa u^3_0}=\sum_k\frac{1}{4\prod_{j\neq k}(u^k-u^j)}=0$ we get
the fact that the families of surfaces $\Phi=\mathrm{ct},\
\frac{\pa\Phi}{\pa u^2_0}=\mathrm{ct},\ \frac{\pa\Phi}{\pa
u^3_0}=\mathrm{ct}$ form a triply orthogonal system and since
$\Phi=\mathrm{ct}$ are parallel we obtain that the last two
families of surfaces are the two families of developables of the
congruence of normals to the first family.

With $\Del(u):=(u-u^2_0)(u-u^3_0)\prod_j(a_j-u)$ and applying $d$
to $\frac{\pa\Phi}{\pa u^2_0}=\mathrm{ct}$ we get
$0=d\frac{\pa\Phi}{\pa
u^2_0}=\frac{1}{4}\sum_k(-1)^k\ep_k\frac{(u^k-u^3_0)du^k}{\sqrt{\Del(u^k)}}$,
so for $u^2\nearrow u^2_0$ we have $du^2=0$; thus the developables
$\frac{\pa\Phi}{\pa u^2_0}=\mathrm{ct}$ are tangent to
$\{u^2=u^2_0\}$ and similarly the developables $\frac{\pa\Phi}{\pa
u^3_0}=\mathrm{ct}$ are tangent to $\{u^3=u^3_0\}$.

The developables $\frac{\pa\Phi}{\pa u^2_0}=\mathrm{ct}$ being
tangent to $\{u^2=u^2_0\}$, have their lines of striction on
$\{u^3=u^3_0\}$ which are geodesics (this is true for any normal
congruence); thus we get Jacobi's equations of geodesics on
$\{u^3=u^3_0\}$ with tangents tangent to $\{u^2=u^2_0\}$
(depending on two constants $u^2_0,\ c$):
\begin{eqnarray}\label{eq:geodell}
\ep_1\int\frac{(u^1-u^3_0)du^1}{\sqrt{\Del(u^1)}}
-\ep_2\int\frac{(u^2-u^3_0)du^2}{\sqrt{\Del(u^2)}}=c.
\end{eqnarray}
Note also that the $1$-dimensional family of geodesics
corresponding to the same value of the constant $u^2_0$ (which can
have here any value between $a_1$ and $a_3$) are tangent to the
same lines of curvature on the ellipsoid $\{u^3=u^3_0\}$. For
example for $u^2=u^2_0$ we have $du^2=0$, so the geodesics are
tangent to the two lines of curvature $(u^2,u^3)=(u^2_0,u^3_0)$ on
the ellipsoid $\{u^3=u^3_0\}$ (the geodesic bounces between them);
if one such geodesic is closed, then all are closed and of same
length.

It is this result of Jacobi's from 1838 (published in 1839) that
opened the whole area of research on confocal quadrics and
(hyper-)elliptic integrals of the XIX$^{\mathrm{th}}$ century
mentioned in the preface of Salmon \cite{S1}.

Finally we are able now to discuss the differential equation of a
line tangent to both $\{u^2=u^2_0\}$ and $\{u^3=u^3_0\}$ and its
linear element in elliptic coordinates: clearly such a line is the
intersection of two developables $\frac{\pa\Phi}{\pa
u^2_0}=\mathrm{ct},\ \frac{\pa\Phi}{\pa u^3_0}=\mathrm{ct}$ and
its arc-length can be taken to be $s:=\Phi$.

We thus have
$0=\sum_k(-1)^k\ep_k\frac{(u^k-u^3_0)du^k}{\sqrt{\Del(u^k)}}
=\sum_k(-1)^k\ep_k\frac{(u^k-u^2_0)du^k}{\sqrt{\Del(u^k)}},\\
-2ds=\sum_k(-1)^k\ep_k\frac{(u^k-u^2_0)(u^k-u^3_0)du^k}{\sqrt{\Del(u^k)}}\Rightarrow
0=\sum_k\frac{(-1)^k\ep_kdu^k}{\sqrt{\Del(u^k)}}=
\sum_ku^k\frac{(-1)^k\ep_kdu^k}{\sqrt{\Del(u^k)}},\\ -2ds
=\sum_k(u^k)^2\frac{(-1)^k\ep_kdu^k}{\sqrt{\Del(u^k)}}$, from
where we get
\begin{eqnarray}\label{eq:linel}
\frac{\ep_1du^1}{(u^2-u^3)\sqrt{\Del(u^1)}}=\frac{\ep_2du^2}{(u^1-u^3)\sqrt{\Del(u^2)}}
=\frac{\ep_3du^3}{(u^1-u^2)\sqrt{\Del(u^3)}}=\frac{2ds}{(u^1-u^2)(u^1-u^3)(u^2-u^3)}.\nonumber\\
\end{eqnarray}
In particular for $u^3=u^3_0$ we have $du^3=0$ and we obtain the
differential equation of geodesics on $\{u^3=u^3_0\}$ with
tangents tangent to $\{u^2=u^2_0\}$ and for
$(u^2,u^3)=(u^2_0,u^3_0)$ we have $du^2=du^3=0$ and we obtain the
differential equation of lines of curvature common to both
$\{u^2=u^2_0\}$ and $\{u^3=u^3_0\}$.

So far we have not paid attention to the signs appearing in the
formulae, since we have not made much use of them. However, in
order to keep a precise accounting of the elliptic coordinates
along threads we must now keep a good accounting of the signs:  as
$s$ increases $\ep_kdu^k$ are always positive, so $\ep_k$ changes
when the variable $u^k$ changes monotonicity.

The four choices of sign $\ep_1\ep_2,\ \ep_1\ep_3$ in the
differential equation (\ref{eq:linel}) of the common tangent
correspond at the level of the geometric picture to the
intersection of tangent cones from a point in space to the two
confocal quadrics being four lines; these lines are obtained by
reflecting one of them in the tangent planes of the three confocal
quadrics that pass through that point and if for example we want
to pass from one of them to the one symmetric with respect to the
tangent plane to the $\{u^j=\mathrm{ct}\}$ quadric, then we change
the sign of $du^j$. These simple remarks allow the generalization
of Chasles's result.

\subsection{Darboux's generalization of Chasles's result}\noindent

We have

\begin{theorem} (Darboux)

Consider a ray of light tangent to two given confocal quadrics
$\{u^j=\mathrm{ct}\},\ \{u^k=\mathrm{ct}\},\ j\neq k$ and which
reflects consecutively in $n$ other quadrics of the confocal
family. If after a certain number of reflections the ray of light
returns to its original position, then this property and the
perimeter of the obtained polygon is independent of the original
position of the ray of light.

\end{theorem}

\begin{proof}

For simplicity we assume that the two initial confocal quadrics
are the hyperboloid with one sheet $\{u^2=u^2_0\}$ and the
ellipsoid $\{u^3=u^3_0\}$, further that $n=1$ and that this
quadric (on which the vertices of the polygon are situated) is an
ellipsoid $\{u^3=u^3_1\}$; thus $u^3_0>u^3_1$, $\{u^3=u^3_0\}$
lies inside $\{u^3=u^3_1\}$, the points of tangency with
$\{u^3=u^3_0\}$ ($\{u^2=u^2_0\}$) are situated on (and possibly
outside) the actual light trajectory segments and on the actual
segments we have $a_1>u^1>a_2>u^2_0>u^2>a_3>u^3_0>u^3>u^3_1$.

Consider the polygon $\mathcal{P}$ with vertices $A_1,A_2,...,A_m$
situated on $\{u^3=u^3_1\}$; since
$\\\sum_k\int\frac{(-1)^{k+1}\ep_kdu^k}{\sqrt{\Del(u^k)}}=$ct we
have
$\sum_k\int_{\mathcal{P}}\frac{(-1)^{k+1}\ep_kdu^k}{\sqrt{\Del(u^k)}}=0$.
As we move from $A_1$ towards $A_2$ $u^3$ increases from $u^3_1$
to $u^3_0$ until we reach the point of tangency with
$\{u^3=u^3_0\}$, then it decreases to $u^3_1$; thus $du^3$ changes
sign from $+$ to $-$ at $\{u^3=u^3_0\}$ and the third integral on
$A_1A_2$ is $2\int_{u^3_1}^{u^3_0}\frac{du^3}{\sqrt{\Del(u^3)}}$;
on the whole polygon $\mathcal{P}$ it is
$2m\int_{u^3_1}^{u^3_0}\frac{du^3}{\sqrt{\Del(u^3)}}$.

For the first two integrals we don't have any change in the sign
of the variations $du^1,\ du^2$ at the vertices of $\mathcal{P}$
(assumed not to be on the planes of coordinates and on
$\{u^2=u^2_0\}$), so they vary freely within their domain of
variation.

Thus we obtain
$\int_{\mathcal{P}}\frac{\ep_1du^1}{\sqrt{\Del(u^1)}}=
2n\int_{a_2}^{a_1}\frac{du^1}{\sqrt{\Del(u^1)}},\ n$ being the
number of times the value $a_2$ is taken by $u^1$, that is the
number of points of intersections of $\mathcal{P}$ with the plane
$\{x^2=0\}$, and thus an even number (for each such occurrence we
double the multiplicity since the value $u^1=a^2$ is taken from
both sides of the plane; note also that $\mathcal{P}$ cuts the
planes $\{x^1=0\},\ \{x^2=0\}$ alternatively).

Similarly $-\int_{\mathcal{P}}\frac{\ep_2du^2}{\sqrt{\Del(u^2)}}=
-2n'\int_{a_3}^{u^2_0}\frac{du^2}{\sqrt{\Del(u^2)}},\ n'$ being
again an even number ($\mathcal{P}$ touches the hyperboloid with
one sheet $\{u^2=u^2_0\}$ and cuts the plane $\{x^3=0\}$
alternatively, since between any two passes through the plane
$\{x^3=0\}$ the other extreme value of $u^2$ is achieved smoothly
and thus at a point of tangency with a hyperboloid with one sheet
which is forced to be $\{u^2=u^2_0\}$).

Thus we have
\begin{eqnarray}\label{eq:darb}
n\int_{a_2}^{a_1}\frac{du^1}{\sqrt{\Del(u^1)}}
-n'\int_{a_3}^{u^2_0}\frac{du^2}{\sqrt{\Del(u^2)}}+
m\int_{u^3_1}^{u^3_0}\frac{du^3}{\sqrt{\Del(u^3)}}=0,\nonumber\\
n\int_{a_2}^{a_1}\frac{u^1du^1}{\sqrt{\Del(u^1)}}
-n'\int_{a_3}^{u^2_0}\frac{u^2du^2}{\sqrt{\Del(u^2)}}+
m\int_{u^3_1}^{u^3_0}\frac{u^3du^3}{\sqrt{\Del(u^3)}}=0,\nonumber\\
2n\int_{a_2}^{a_1}\frac{(u^1)^2du^1}{\sqrt{\Del(u^1)}}
-2n'\int_{a_3}^{u^2_0}\frac{(u^2)^2du^2}{\sqrt{\Del(u^2)}}+
2m\int_{u^3_1}^{u^3_0}\frac{(u^3)^2du^3}{\sqrt{\Del(u^3)}}=\mathrm{perimeter}(\mathcal{P})
\end{eqnarray}
the last two relations being obtained by similar computations. The
first two relations impose certain rationality conditions on
hyper-elliptic integrals involving $u^2_0,\ u^3_0,\ u^3_1$ (so in
general such configurations do not exist); the last relation shows
that the perimeter of $\mathcal{P}$ is the same for all such
polygons $\mathcal{P}$.

All that remains to check is that the existence of such polygons
$\mathcal{P}$ is an open condition.

Consider $B_1$ infinitesimally close to $A_1$ on the ellipsoid
$\{u^3=u^3_1\}$; through it we draw a line tangent to both
$\{u^2=u^2_0\}$ and $\{u^3=u^3_0\}$ and infinitesimally close to
$A_1A_2$ to obtain $B_2\in\{u^3=u^3_1\}$ infinitesimally close to
$A_2$; etc. Thus we obtain the polygonal line $\mathcal{P}':\
B_1,B_2,...,B_m,B'_1$ and we want $B'_1=B_1$.

The computations from the first two relations of (\ref{eq:darb})
still hold, except for a small deficit in the variations of the
first two elliptic coordinates near the coordinates $(u^1,u^2)$ of
$B_1$ and $(u'^1,u'^2)$ of $B'_1$:
$$\frac{du^1}{\sqrt{\Del(u^1)}}-\frac{du^2}{\sqrt{\Del(u^2)}}=
\frac{du'^1}{\sqrt{\Del(u'^1)}}-\frac{du'^2}{\sqrt{\Del(u'^2)}},\
\frac{u^1du^1}{\sqrt{\Del(u^1)}}
-\frac{u^2du^2}{\sqrt{\Del(u^2)}}=
\frac{u'^1du'^1}{\sqrt{\Del(u'^1)}}
-\frac{u'^2du'^2}{\sqrt{\Del(u'^2)}}.$$ This can be interpreted as
a non-degenerate differential system in $(u'^1,u'^2)$; since for
$B_1=A_1$ it has the initial condition $(u'^1,u'^2)=(u^1,u^2)$, we
conclude that the obvious solution $(u'^1,u'^2)=(u^1,u^2)$ is
unique.

\end{proof}

\begin{remark}

Note that if we allow the ellipsoid $\{u^3=u^3_0\}$ in Darboux's
result to become degenerated (that is $u^3_0=a_3$) and one of the
vertices of the polygon is situated in the $\{x^3=0\}$ plane, then
we get Chasles's result; in this vein Darboux's generalization of
Chasles's result is a statement about closed geodesics on a
degenerated $3$-dimensional ellipsoid in $\mathbb{R}^4$.

\end{remark}

\subsection{Staude's generalization of Graves's result}\noindent

Note that the computations in the last part of (\ref{eq:darb})
referring to arc-length are valid also for line of curvature
$(u^2,u^3)=(u^2_0,u^3_0)$ segments and for geodesics on
$\{u^3=u^3_0\}$ with tangents tangent to $\{u^2=u^2_0\}$ segments.

Staude's thread construction of confocal ellipsoids roughly states

\begin{theorem} (Staude)

If a thread $\mathcal{P}$ passed around an ellipsoid
$\{u^3=u^3_0\}$ and outside the hyperboloid with one sheet
$\{u^2=u^2_0\}$ be stretched with a pen at a point $P$ such that
it touches once both visible parts of the hyperboloid with one
sheet $\{u^2=u^2_0\}$ (thus the thread consists either of two
rectilinear segments from $P$ tangent to both $\{u^2=u^2_0\}$ and
$\{u^3=u^3_0\}$, two geodesic segments on $\{u^3=u^3_0\}$ with
tangents tangent to $\{u^2=u^2_0\}$ and in continuation of the two
rectilinear segments, two line of curvature
$(u^2,u^3)=(u^2_0,u^3_0)$ segments in continuation of the previous
two geodesic segments and another geodesic segment on
$\{u^3=u^3_0\}$ with tangents tangent to $\{u^2=u^2_0\}$ and in
continuation of and joining the two line of curvature
$(u^2,u^3)=(u^2_0,u^3_0)$ segments, or, if a rectilinear segment
of the thread touches the hyperboloid with one sheet
$\{u^2=u^2_0\}$ before the ellipsoid $\{u^3=u^3_0\}$, then it is
not required to touch the corresponding line of curvature branch),
then the point $P$ will move on an ellipsoid $\{u^3=u^3_1\}$
confocal to and outside $\{u^3=u^3_0\}$ and $\{u^2=u^2_0\}$.

\end{theorem}

\begin{proof}

Note that with $n,\ n',\ m$ being the same as in Darboux's result
we have $n=n'=2,\ m=1$, since the thread cuts alternatively the
planes $\{x^1=0\},\ \{x^2=0\}$ twice, touches the hyperboloid with
one sheet $\{u^2=u^2_0\}$ and cuts the plane $\{x^3=0\}$
alternatively twice (here we use the fact that the thread touches
once each visible part of the hyperboloid with one sheet
$\{u^2=u^2_0\}$) and touches the ellipsoid $\{u^3=u^3_1\}$ on
which $P$ is situated only once (at $P$).

For the length of the thread $\mathcal{P}$ we get (similarly to
the last relation of (\ref{eq:darb}))

\begin{eqnarray}\label{eq:staud}
4\int_{a_2}^{a_1}\frac{(u^1-u^2_0)(u^1-u^3_0)du^1}{\sqrt{\Del(u^1)}}
-4\int_{a_3}^{u^2_0}\frac{(u^2-u^2_0)(u^2-u^3_0)du^2}{\sqrt{\Del(u^2)}}+
2\int_{u^3_1}^{u^3_0}\frac{(u^3-u^2_0)(u^3-u^3_0)du^3}{\sqrt{\Del(u^3)}},\nonumber\\
\end{eqnarray}
so it depends only on the third elliptic coordinate $u^3_1$ of
$P$.

\end{proof}

\begin{remark}

For $u^2_0=a_2$ or $u^3_0=a_3$ (or both) the line of curvature
segments become geodesic segments and we have a genuine
variational problem. For $(u^2_0,u^3_0)=(a_2,a_3)$ the thread
configuration as it appears in Hilbert-Cohn-Vossen (\cite{HCV},\S
4) and Salmon (\cite{S1},\S 421a,b) changes in the sense that one
of the rectilinear segments is not continued until it reaches
$\{u^3_0=a_3\}$ after touching $\{u^2_0=a_2\}$, but it is broken
to reach the focus of the corresponding branch of the focal
hyperbola (there is no well defined normal for the degenerated
hyperboloid with one sheet $\{u^2_0=a_2\}$ at the points of its
singular boundary-hyperbola, so this construction makes sense from
a geometric point of view). The two points of view are equivalent
and are explained by the thread construction of an ellipse with
the thread passing through its foci being able to be replaced with
the same construction, but the foci being replaced with any point
of the corresponding branch of the hyperbola which is orthogonal
focal curve of the ellipse.

\end{remark}

\section{Thread configurations for ellipsoids}

Just as Darboux, in order to simplify the presentation we assume
that the vertices of the thread configuration move on the same
ellipsoid confocal to the given one.

Since the geodesic segments in the thread are in continuation of
rectilinear segments in the thread, in order to allow variations
of vertices on their corresponding confocal ellipsoids we invert
the point of view: rectilinear segments are in continuation of
geodesic segments. Thus we need to consider the intersection of a
forward and a backward part of tangent surfaces of two geodesic
segments on $\{u^3=u^3_0\}$ with tangents tangent to
$\{u^2=u^2_0\}$: this intersection is a curve situated on a
confocal ellipsoid and the length of a thread fixed at two points
on the two geodesic segments and stretched with a pen situated on
this curve is constant.

Conversely, if the forward part of the tangent surface  of a
geodesic segment on $\{u^3=u^3_0\}$ with tangents tangent to
$\{u^2=u^2_0\}$ is reflected in a confocal ellipsoid
$\{u^3=u^3_1<u^3_0\}$, then we obtain the backward part of the
tangent surface of a geodesic segment on $\{u^3=u^3_0\}$ with
tangents tangent to $\{u^2=u^2_0\}$.

Thus we are led to the basic local result which allows thread
configurations for ellipsoids and which appears (implicitly) in
both Darboux's and Staude's results:

\begin{proposition} (Staude's vertex configuration)

The forward and backward parts of the tangent surfaces of two
geodesic segments on $\{u^3=u^3_0\}$ with tangents tangent to
$\{u^2=u^2_0\}$ meet along a confocal ellipsoid
$\{u^3=u^3_1>u^3_0\}$; moreover the length of the thread measured
from a fixed point on the first geodesic to a fixed point on the
second geodesic and via the vertex on the confocal ellipsoid
$\{u^3=u^3_1\}$ is constant.

\end{proposition}

\begin{proof}

We have such tangent surfaces $\frac{\pa\Phi}{\pa
u^2_0}=\frac{1}{4}\sum_k(-1)^k\int\frac{(u^k-u^3_0)\ep_kdu^k}{\sqrt{\Del(u^k)}}=\mathrm{ct},\
\ep_kdu^k>0$; $\ep_1$ changes sign at the intersection with the
planes of coordinates $\{x^1=0\},\ \{x^2=0\}$, $\ep_2$ changes
sign at the intersection with the plane of coordinate $\{x^3=0\}$
and at the points of tangency with $\{u^2=u^2_0\}$ and $\ep_3$
changes sign from $-$ on the backward part of the tangent surface
to $+$ on the forward part of the tangent surface along the
geodesic on $\{u^3=u^3_0\}$. If we reflect the forward part of the
tangent surface along the ellipsoid $\{u^3=u^3_1\}$, then $\ep_3$
changes sign from $+$ to $-$ and we get the backward part of a
similar tangent surface. The length of the thread measured from a
fixed point on the first geodesic to a fixed point on the second
one is constant since $u^1,\ u^2$ vary between given values and
$u^3$ varies twice between $u^3_1$ and $u^3_0$.

\end{proof}

Since the computations in Darboux's result are mostly valid also
for geodesics on $\{u^3=u^3_0\}$ with tangents tangent to
$\{u^2=u^2_0\}$, they remain valid if we allow such geodesic
segments in the thread, in which case (\ref{eq:darb}) is replaced
with
\begin{eqnarray}\label{eq:darb1}
n\int_{a_2}^{a_1}\frac{(u^1-u^3_0)du^1}{\sqrt{\Del(u^1)}}
-n'\int_{a_3}^{u^2_0}\frac{(u^2-u^3_0)du^2}{\sqrt{\Del(u^2)}}+
m\int_{u^3_1}^{u^3_0}\frac{(u^3-u^3_0)du^3}{\sqrt{\Del(u^3)}}=0,\nonumber\\
2n\int_{a_2}^{a_1}\frac{u^1(u^1-u^3_0)du^1}{\sqrt{\Del(u^1)}}
-2n'\int_{a_3}^{u^2_0}\frac{u^2(u^2-u^3_0)du^2}{\sqrt{\Del(u^2)}}+
2m\int_{u^3_1}^{u^3_0}\frac{u^3(u^3-u^3_0)du^3}{\sqrt{\Del(u^3)}}
=\mathrm{length}(\mathcal{P}).\nonumber\\
\end{eqnarray}
Note that (\ref{eq:darb}) implies (\ref{eq:darb1}): as the sides
$A_1A_2,\ A_2A_3,...,A_mA_1$ move on their subjacent tangent
surfaces of geodesic segments on $\{u^3=u^3_0\}$ the polygon
$\mathcal{P}$ closes and its length remains constant. Conversely,
if the polygon fails to close, but $A_m,A_1$ still can be joined
by a thread formed by two rectilinear segments and a geodesic
segment on $\{u^3=u^3_0\}$ with tangents tangent to
$\{u^2=u^2_0\}$ in continuation of the two rectilinear segments
and joining them, then only the rationality condition from
(\ref{eq:darb1}) remains valid among the rationality conditions
(\ref{eq:darb}). The length of the geodesic segment in the thread
$A_mA_1$ can be liberally distributed among the segments
$A_1A_2,...A_{m-1}A_m$ (thus making them pieces of thread); if we
allow consecutive rectilinear segments in the thread $\mathcal{P}$
to be joined by geodesic segments but with cusps at both ends,
then we can prescribe any desired length to all but one of the
geodesic segments in the thread $\mathcal{P}$ (of course the
integrals along the geodesic segments with cusps at ends have to
be subtracted).

Summing up, the rationality condition of (\ref{eq:darb1}) is the
condition that after several iteration of the local result (the
backward part of the tangent surface obtained by reflection in the
ellipsoid $\{u^3=u^3_1\}$ of the forward part of another tangent
surface being continued to obtain the forward part of the tangent
surface and thus allowing the iteration) the process closes up.
Since we have only an $1$-dimensional family of geodesic segments
under consideration, only a rationality condition is required as
closing condition, so the first equation of (\ref{eq:darb1})
provides also a sufficient condition for an open set of such
thread configurations.

In particular if we have no vertices ($m=0$), then the first
equation of (\ref{eq:darb1}) becomes the rationality condition for
the existence of closed geodesics on $\{u^3=u^3_0\}$ with tangents
tangent to $\{u^2=u^2_0\}$ and the last equation gives the length
of such closed geodesics.

Thus Darboux's result (having no geodesic segments) and closed
geodesics (having only geodesic segments) are the two extremes of
thread configurations formed by rectilinear and geodesic segments.

For Staude's result there is still the question of the thread
closing up at $P$; thus we need to consider the question for
thread without line of curvature segments; from the first equation
of (\ref{eq:darb1}) this is equivalent to
\begin{eqnarray}\label{eq:halfgeod}
\int_{a_2}^{a_1}\frac{(u^1-u^3_0)du^1}{\sqrt{\Del(u^1)}}
-\int_{a_3}^{u^2_0}\frac{(u^2-u^3_0)du^2}{\sqrt{\Del(u^2)}}(=
-\frac{1}{2}\int_{u^3_1}^{u^3_0}\frac{(u^3-u^3_0)du^3}{\sqrt{\Del(u^3)}})>0.
\end{eqnarray}
We are interested in this relation since it has to do with the
variation of elliptic coordinates on geodesics: consider a
geodesic segment on $\{u^3=u^3_0\}$ with tangents tangent to
$\{u^2=u^2_0\}$ touching once each branch of the line of curvature
$(u^2,u^3)=(u^2_0,u^3_0)$ (the geodesic segment begins and ends at
such touching points). From (\ref{eq:darb1}) this happens in more
than half a turn around the $e_3$-axis (that is the plane through
the $e_3$-axis and passing through one of the touching points cuts
again the considered geodesic segment) iff we have
(\ref{eq:halfgeod}) (note that for $u^2_0=a_2$ equality is
obtained in (\ref{eq:halfgeod}) since (\ref{eq:geodell}) forces
geodesics passing through an umbilic of an ellipsoid to also pass
through its opposite one). If we extend this geodesic segment at
each end with similar geodesic segments, then the forward
(backward) half of the tangent surfaces of these new two geodesic
segments will cut along a curve on a confocal ellipsoid; for $P$
inside this ellipsoid Staude's thread construction is not possible
(of course one can correct this by allowing a longer line of
curvature segment, in which case we have $n=4$ or by allowing an
ideal construction with a smaller line of curvature segment with
cusps at both ends and thus its length must be subtracted); for
$P$ outside (on) this ellipsoid the thread construction is
possible and does (not) contain line of curvature segments.

If in (\ref{eq:halfgeod}) we have the opposite inequality (for
$u^3_0$ close to $a_3$ such a geodesic segment looks close enough
to the two tangent segments from Graves's result), then the story
changes completely: Staude's thread construction is possible for
all confocal ellipsoids and always it contains line of curvature
segments.

If the rationality condition from (\ref{eq:darb1}) is not
satisfied, then one can make again liberal use of line of
curvature $(u^2,u^3)=(u^2_0,u^3_0)$ segments in the thread
(including with cusps at both ends) when the geodesic segments of
the thread touch such lines of curvature; all but one such line of
curvature segments can be arbitrarily prescribed and the length of
the thread $\mathcal{P}$ is
\begin{eqnarray}\label{eq:staud1}
2n\int_{a_2}^{a_1}\frac{(u^1-u^2_0)(u^1-u^3_0)du^1}{\sqrt{\Del(u^1)}}
-2n'\int_{a_3}^{u^2_0}\frac{(u^2-u^2_0)(u^2-u^3_0)du^2}{\sqrt{\Del(u^2)}}+
m\int_{u^3_1}^{u^3_0}\frac{(u^3-u^2_0)(u^3-u^3_0)du^3}{\sqrt{\Del(u^3)}}.\nonumber\\
\end{eqnarray}

\section{Darboux's Theorem in general setting}

Since quadrics with distinct non-zero eigenvalues of the quadratic
part defining the quadric form an open dense set in the set of all
quadrics, by a continuity argument one can infer that Darboux's
generalization of Chasles's result is valid for all quadrics.

However, while elliptic coordinates a-priori are a must for
geodesics and lines of curvature on quadrics, they should not be
necessary for straight line segments and should be replaced by
purely algebraic computations; locally the Ivory affinity provides
such a venue and Darboux's constant perimeter property concerning
moving polygons circumscribed to a given set of $n$ and inscribed
in arbitrarily many $n$-dimensional confocal quadrics provides the
needed global arguments.

\subsection{An Ivory affinity approach for the vertex
configuration}\noindent

\noindent
We have now

\begin{proposition} (Ivory affinity approach for the vertex configuration)

{\it Given three points $x_0^0,x_0^1,x_0^2\in x_0$ and by the
Ivory affinity the corresponding points $x_z^0,x_z^1,x_z^2\in
x_z$, then $V_1^0,\ V_2^0$ reflect in $x_z$ at $x_z^0$ iff
$V_0^1,\ V_0^2$ reflect in $x_0$ at $x_0^0$; further in this case,
by the preservation of the TC under the Ivory affinity, we
conclude that $V_1^0$ is tangent to $x_0$ at $x_0^1$ iff $V_2^0$
is tangent to $x_0$ at $x_0^2$ and further in this case
$x_z^1,x_0^0,x_z^2$ are co-linear. Also $V_0^1,\ V_1^0$ are
tangent to the same set of quadrics confocal to the given one
$x_0$.}

\end{proposition}

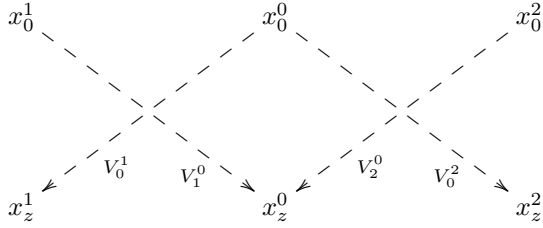
\begin{center}
$\xymatrix@!0{x_0^1\ar@{-->}[dddrrrr]_>>>>>>>>{V_1^0}&&&&
x_0^0\ar@{-->}[dddllll]^>>>>>>>>>>{V_0^1}\ar@{-->}[dddrrrr]_>>>>>>>>{V_0^2}&&&&
x_0^2\\  \\  \\
x_z^1&&&& x_z^0\ar@{<--}[uuurrrr]_<<<<<<<<<<{V_2^0}&&&& x_z^2}$
\end{center}

\begin{remark}

Note that excepting totally real cases one loses the orientation
of the real numbers, so from a complex point of view refractions
must also be considered as reflections.

\end{remark}

\begin{proof}

Note $\sqrt{R_z}V_1^0=-V_0^1-z\hat N_0^0$ (here we use $\hat
N_0^0=Ax_0^0+B,\ (I_{n+1}+\sqrt{R_z})C(z)+zB=0$).

We have
$(dx_z^0)^T(\frac{V_1^0}{|V_1^0|}\pm\frac{V_2^0}{|V_2^0|})=(dx_0^0)^T\sqrt{R_z}
(\frac{V_1^0}{|V_1^0|}\pm\frac{V_2^0}{|V_2^0|})
=-(dx_0^0)^T(\frac{V_0^1}{|V_0^1|}\pm\frac{V_0^2}{|V_0^2|})$.

Thus $V_1^0,V_2^0$ reflect in $x_z$ at $x_z^0$ iff $V_0^1,V_0^2$
reflect in $x_0$ at $x_0^0$ and the tangency part follows from
this.

The vector $V_0^1$ is tangent to $x_{z'}$ iff the quadratic
equation in $t$
$0=Q_{z'}(x_0^0+tV_0^1)=(V_0^1)^TAR_{z'}^{-1}V_0^1t^2\\+2(V_0^1)^TR_{z'}^{-1}\hat
N_0^0t+z'(\hat N_0^0)^TR_{z'}^{-1}\hat N_0^0$ has double root and
we would like the same to hold for $Q_{z'}(x_0^1+tV_1^0)=0$, that
is we want the discriminant $\Del:=[(V_0^1)^TR_{z'}^{-1}\hat
N_0^0]^2-z'(V_0^1)^TAR_{z'}^{-1}V_0^1(\hat N_0^0)^TR_{z'}^{-1}\hat
N_0^0$ to be symmetric in $x_0^0,\ x_0^1$. This is so because
$\sqrt{R_z}V_0^1=-V_1^0-z\hat N_0^1$, so $z^2\Del=
[(V_0^1)^TR_{z'}^{-1}(V_0^1+\sqrt{R_z}V_1^0)]^2
+(V_0^1)^T(I_{n+1}-R_{z'}^{-1})V_0^1(V_0^1+\sqrt{R_z}V_1^0)^TR_{z'}^{-1}(V_0^1+\sqrt{R_z}V_1^0)
=[(V_0^1)^TR_{z'}^{-1}\sqrt{R_z}V_1^0]^2+|V_0^1|^2[2(V_0^1)^TR_{z'}^{-1}\sqrt{R_z}V_1^0
+(1-\frac{z}{z'})[(V_0^1)^TR_{z'}^{-1}V_0^1
+(V_1^0)^TR_{z'}^{-1}V_1^0]]
-(1-\frac{z}{z'})(V_0^1)^TR_{z'}^{-1}V_0^1(V_1^0)^TR_{z'}^{-1}V_1^0\\
+\frac{z}{z'}|V_0^1|^2|V_1^0|^2$ and from the Ivory Theorem we
have $|V_0^1|^2=|V_1^0|^2$.

\end{proof}

Now we are able to derive the vertex configuration of Darboux's
Theorem: assume that $V_1^0$ is tangent to $x_0$ at $x_1^0$; by
reflection in $x_z$ at $x_z^0$ we get $V_2^0$ tangent to $x_0$ at
$x_0^2$. If $V_1^0$ is tangent to $x_{z'}$, then so is $V_0^1$;
since $V_0^1$ and $V_0^2$ are co-linear, the same statement holds
for $V_0^2$, so we finally conclude that $V_2^0$ is tangent to
$x_{z'}$.

\subsection{Darboux's Theorem in general setting}

\begin{theorem} (Darboux in general setting)

Consider a ray of light tangent to $n$ given confocal quadrics in
$\mathbf{C}^{n+1}$ and which reflects consecutively in $p$ other
given quadrics of the confocal family. If after a minimal number
$p$ of reflections the ray of light returns to its original
position, then this property and the perimeter of the obtained
polygon is independent of the original position of the ray of
light.

\end{theorem}

\begin{proof}

Consider the polygonal line $\mathcal{P}$ with vertices
$x_{z_j}^j\in Q_{z_j},\ j=0,1,2,...,p$ and which is obtained as
follows: pick arbitrarily $x_{z_0}^0\in Q_{z_0}$, then
$x_{z_1}^1\in Q_{z_1}$ such that $x_{z_1}^1-x_{z_0}^0$ is tangent
to the n given confocal quadrics  $Q_{z'_1},...,\ Q_{z_n'}$ at
points $y_{z'_1}^1$ (according to Chasles-Jacobi's result and an
observation of Darboux's there are choices of the direction
$x_{z_1}^1-x_{z_0}^0$ given by the intersections of the tangent
cones to $Q_{z'_1},...,\ Q_{z_n'}$ through $x_{z_0}^0$; note also
that for each choice of direction there are in general two choices
of $x_{z_1}^1$); now choose $x_{z_2}^2\in Q_{z_2}$ such that
$x_{z_1}^1-x_{z_0}^0,\ x_{z_2}^2-x_{z_1}^1$ reflect in $Q_{z_1}$
at $x_{z_1}$ (note that for $z_2\neq z_0$ there are in general two
choices of $x_{z_2}^2$ such that $x_{z_2}^2\neq x_{z_0}^0$, since
$x_{z_2}^2-x_{z_1}^1$ will not be tangent in general to
$Q_{z_2}$), etc.

The polygonal line $\mathcal{P}$ closes to a Darboux $p$-polygon
iff $z_0=z_p,\ x_{z_0}^0=x_{z_p}^p,\ x_{z_p}^p-x_{z_{p-1}}^{p-1},\
x_{z_1}^1-x_{z_0}^0$ reflect in $Q_{z_0}$ at $x_{z_0}^0$ and $p$
is minimal (a picture with $p=3$ should clarify the issues under
discussion).

Now by Chasles-Jacobi's result and {\bf Proposition 6.1} in
general each segment $x_{z_j}^j-x_{z_{j-1}}^{j-1},\ j=1,...,p$ of
the polygonal line $\mathcal{P}$ is tangent to same $n$ quadrics
$Q_{z'_1},...,\ Q_{z_n'}$ at points $y_{z'_k}^j,\ k=1,...,n$.

The condition that $\mathcal{P}$ closes to a Darboux $p$-polygon
is not just the condition that $x_{z_0}^0-x_{z_{p-1}}^{p-1}$ is
tangent to the $n$ quadrics $Q_{z'_1},...,\ Q_{z_n'}$ at points
$y_{z'_k}^0,\ k=1,...,n$: by Chasles-Jacobi's such a line is
uniquelly determined by being tangent to $n$ quadrics, having a
passing point and a direction chosen among the $(n+1)$ directions
(reflections in the principal spaces passing through that point;
they are decided by the differentials of the $(n+1)$ points moving
on the $(n+1)$ quadrics passing through that point); knowing this
direction among the $(n+1)$ choices will force the reflection
property needed to close the Darboux $p$-polygon.

The data about the polygonal line $\mathcal{P}$ we have is
\begin{eqnarray}\label{eq:tare}
(dx_{z_j}^j)^T
(\ep_j\frac{x_{z_j}^j-x_{z_{j-1}}^{j-1}}{|x_{z_j}^j-x_{z_{j-1}}^{j-1}|}
-\ep_{j+1}\frac{x_{z_{j+1}}^{j+1}-x_{z_j}^j}{|x_{z_{j+1}}^{j+1}-x_{z_j}^j|})=0,\
\ep_j:=\pm 1,\ j=1,...,p-1
\end{eqnarray}
and we get $\mathcal{P}$ closing to a Darboux $p$-polygon  iff
\begin{eqnarray}\label{eq:reta}
\ \ \ \ \ z_0=z_p,\ x_{z_0}^0=x_{z_p}^p,\ \ep_0=\ep_p:=\pm 1,\
p\ \mathrm{minimal},\nonumber\\
(dx_{z_p}^p)^T
(\ep_p\frac{x_{z_p}^p-x_{z_{p-1}}^{p-1}}{|x_{z_p}^p-x_{z_{p-1}}^{p-1}|}-
\ep_1\frac{x_{z_1}^1-x_{z_0}^0}{|x_{z_1}^1-x_{z_0}^0|})=0.
\end{eqnarray}
Similarly to Darboux's original proof, what we need now to prove
is that (\ref{eq:tare}), (\ref{eq:reta}) and the perimeter
\begin{eqnarray}\label{eq:peri}
\sum_{j=1}^p\ep_j|x_{z_j}^j-x_{z_{j-1}}^{j-1}|=c_0
\end{eqnarray}
of $\mathcal{P}$ are independent of the choice of $x_{z_0}'^0\in
Q_{z_0},\ x_{z_1}'^1\in Q_{z_1}$ respectively near and instead of
$x_{z_0}^0\in Q_{z_0},\ x_{z_1}^1\in Q_{z_1}$ and such that
$z_1,...,z_p,\ \ep_1,...,\ep_p,\ z'_1,...,z'_n$ remain fixed.

If we apply $d$ to (\ref{eq:peri}), then we get the sum of
(\ref{eq:tare}) and (\ref{eq:reta}):
$d\sum_{j=1}^p\ep_j|x_{z_j}^j-x_{z_{j-1}}^{j-1}|=\\\sum\circlearrowleft_{j=1}^p(dx_{z_j}^j)^T
(\ep_j\frac{x_{z_j}^j-x_{z_{j-1}}^{j-1}}{|x_{z_j}^j-x_{z_{j-1}}^{j-1}|}
-\ep_{j+1}\frac{x_{z_{j+1}}^{j+1}-x_{z_j}^j}{|x_{z_{j+1}}^{j+1}-x_{z_j}^j|})
$. The same behaviour appplies to (\ref{eq:tare})' (which is
obtained by construction and induction on $j=1,...,p-1$) and
(\ref{eq:reta})': together are equivalent to (\ref{eq:peri})'
being constant (by continuity presumably the constant of
(\ref{eq:peri})).

Similarly to Darboux's original proof and by a continuity argument
we need only prove $x_{z_p}'^p=x_{z_0}'^0$ and (\ref{eq:reta})';
using this and (\ref{eq:tare})' we get (\ref{eq:peri})'.

Summing up: given $x_{z_j}^j,\ y_{z_k'}^j,\ j=0,...,p,\ k=1,...,n$
as before (to satisfy including (\ref{eq:tare}), (\ref{eq:reta}),
(\ref{eq:peri})) and $'$ quantities constructed from
$x_{z_0}'^0\in Q_{z_0},\ x_{z_1}'^1\in Q_{z_1}$ respectively near
and instead of $x_{z_0}^0\in Q_{z_0},\ x_{z_1}^1\in Q_{z_1}$ and
such that $z_1,...,z_p,\ \ep_1,...,\ep_p,\ z'_1,...,z'_n$ remain
fixed, then all $'$ quantities will be near their corresponding
counterpart quantities and will satisfy, except for $x_{z_p}'^p$,
the relations satisfied by their counterparts. But now by
continuity Chasles-Jacobi's result will be precise enough to
choose from among the choices of $x_{z_p}'^p$ the correct one such
that $x_{z_p}'^p$ will be near $x_{z_p}^p=x_{z_0}^0$.

Now we have the two differential systems
$$(dx_{z_p}^p)^T
(\ep_p\frac{x_{z_p}^p-x_{z_{p-1}}^{p-1}}{|x_{z_p}^p-x_{z_{p-1}}^{p-1}|}-
\ep_1\frac{x_{z_1}^1-x_{z_0}^0}{|x_{z_1}^1-x_{z_0}^0|})=0,$$
$$(dx_{z_p}'^p)^T
(\ep_p\frac{x_{z_p}'^p-x_{z_{p-1}}'^{p-1}}{|x_{z_p}'^p-x_{z_{p-1}}'^{p-1}|}-
\ep_1\frac{x_{z_1}'^1-x_{z_0}'^0}{|x_{z_1}'^1-x_{z_0}'^0|})=0.$$

Since the first system is Pffaf non-exact with unique solution
$x_{z_p}^p=x_{z_0}^0$ near $x_{z_0}^0$ (should it be exact in some
variables, we would have a continuous $(1\le n-1)$-dimensional
family of solutions), by continuity in parameters and nearness of
computations we deduce that also the second system is Pffaf
non-exact with unique solution $x_{z_p}'^p=x_{z_0}'^0$.

\end{proof}

\subsection{An Ivory affinity approach for
newly obtained Darboux $p$-polygons and its iteration}

Take any quadric (for example $Q_{z'_1}$; we can consider
$z'_1:=0$) of the $n$ ones $Q_{z'_k},\ k=1,...,n$. By the Ivory
afffinity as in {\bf Proposition 6.1} one can reverse the r\^{o}le
of the vertices of the polygonal line $\mathcal{P}$ and of the
points of tangency $y_{z'_1}^1,...,y_{z'_1}^p$ with $Q_{z'_1}$,
thus getting a-priori a new Darboux $p$-poligon $\mathcal{\ti{P}}$
with vertices $y_{z_1}^1,...,y_{z_p}^p$ (note that in order to
apply {\bf Proposition 6.1} we need the points
$x_{z_1}^p,x_{z_2}^1...,x_{z_p}^{p-1}$ obtained by the Ivory
affinity), the poligonal segments tangent to $Q_{z'_1}$ at
$x_{z'_1}^1,...,x_{z'_1}^p$ and to $Q_{z'_k}$ at points $\ti
y_{z'_k}^j,\ j=1,...,p,\ k=2,...,n$ different from the original
ones (because they lie on different lines). These Darboux
$p$-polygons are actually closed because the original one
$\mathcal{P}$ is closed  and because we have reflection properties
at all vertices; also they have the same perimeter as that of
$\mathcal{P}$ for $z_1=...=z_p$ (because of the Ivory Theorem on
preservation of lengths of segments between confocal quadrics);
for the remaining general case they have the same perimeter by the
general form of the Darboux Theorem we just proved (note however
that the Darboux polygon $\mathcal{P}$ with the same perimeter may
be different from the original one, as there are many choices in
its construction).

This process can be of course iterated. Does it stop after a
certain number of iterations? If yes, then how many do we actually
get?


\begin{thebibliography}{99}
\def\topset{0pt}
\def\parsep{0pt plus 5pt minus 1pt}
\def\itemsep{-0.5ex}
\small

\bibitem{AF} S. Abenda, Y. Fedorov {\it Closed Geodesics and Billiards on Quadrics related to
elliptic KdV solutions,} preprint
{\href{http://arxiv.org/abs/nlin/0412034} {arXiv:nlin/0412034v1}}
[nlin.SI] 14 Dec 2004.

\bibitem{B} L. Bianchi {\it Lezioni Di Geometria Differenziale,} Vol
{\href{http://visualiseur.bnf.fr/ark:/12148/bpt6k99688w} {{\bf
2}}}, Enrico Spoerri Libraio-Editore, Pisa (1908).

\bibitem{C} J. L. Coolidge {\href{http://www.gutenberg.org/etext/26373}
{\it The Elements of Non-Euclidean Geometry,}} Clarendon Press,
Oxford (1909).

\bibitem{D1} G. Darboux  {\href{http://fr.dleex.com/details/?11194}
{\it Le\c{c}ons Sur La Th\'{e}orie G\'{e}n\'{e}rale Des
Surfaces,}} Vol {\bf 1-4}, Gauthier-Villars, Paris (1894-1917).

\bibitem{D2} G. Darboux {\href{http://fr.dleex.com/details/?11198}
{\it Le\c{c}ons Sur Les Syst\`{e}mes Orthogonaux Et Les
Coordon\'{e}es Curvilignes,}} Gauthier-Villars, Paris (1910).

\bibitem{DR2} V. Dragovi\'{c}, M. Radnovi\'{c} {\it Hyperelliptic Jacobians as Billiard Algebra
of Pencils of Quadrics: Beyond Poncelet Porisms,} preprint
{\href{http://arxiv.org/abs/0710.3656?context=math-ph}
{arXiv:0710.3656v2}} [math.AG] 9 Jul 2008.

\bibitem{F} Y. Fedorov {\it Algebraic Closed Geodesics on a Triaxial Ellipsoid,} preprint
{\href{http://arxiv.org/abs/nlin/0506063} {arXiv:nlin/0506063v2}}
[nlin.SI] 10 Oct 2005.

\bibitem{HCV} D. Hilbert, S. Cohn-Vossen {\it Geometry and Imagination,} Chelsea, New York (1999).

\bibitem{S1} G. Salmon {\it A Treatise On The Analytic Geometry of Three Dimensions,}
 Vol {\bf {\href{http://www.archive.org/details/treatiseonanalyt01salmuoft}
1}-2}, revised by Reginald A. P. Rogers, Longmans, Green and Co.,
New York (1912).

\bibitem{S2} O. Staude {\it Geometrische Deutung der Additionstheoreme der hyperelliptischen
Integrale und Functionen 1. Ordnung im System der confocalen
Fl\"{a}chen 2. Grades,} Math. Ann. XXII,
{\href{http://www.digizeitschriften.de/index.php?id=loader&tx_jkDigiTools_pi1[IDDOC]=300133}
{1-69, 145-176}} (1883).

\bibitem{T} S. Tabachnikov {\href{http://www.math.psu.edu/tabachni/Books/billiardsgeometry.pdf}
{\it Geometry and Billiards}}.

\end{thebibliography}
\end{document}